\documentclass[10pt,a4paper,twocolumn]{IEEEtran}

\IEEEoverridecommandlockouts                              
\usepackage{amsmath}
\usepackage{color}
\usepackage{graphicx}
\usepackage{latexsym}

\usepackage{epsfig} 

\title{On the Equivalence among Three Controllability Problems for a Networked System}

\author{Tong Zhou 
\thanks{This work was supported in part by the NSH under Grant 61573209, 61273126, 65573156  and 51361135705.}
\thanks{T.Zhou is with the Department of Automation and TNList, Tsinghua University, Beijing, 100084,
CHINA. (Tel: 86-10-62797430; Fax: 86-10-62786911; e-mail:
tzhou@mail.tsinghua.edu.cn.)}
}

\begin{document}
\renewcommand{\thefootnote}{\fnsymbol{footnote}}
\maketitle 
\renewcommand{\thepage}{7--\arabic{page}}
\setcounter{page}{1}

\begin{abstract}
A new proof is given for the mathematical equivalence among three $k$-sparse controllability problems of a networked system, which plays key roles in \cite{Olshevsky14} in the establishment of the NP-hardness of the associated minimal controllability problems (MCP). Compared with the available ones, a completely deterministic approach is adopted. Moreover, only primary algebraic operations are utilized in all the derivations. These results enhance the available conclusions about the NP-hardness of a MCP, and can also be directly applied to the computational complexity analysis for a minimal observability problem. In addition, the results of \cite{Olshevsky14} have also been extended to situations in which there are also some other restrictions, such as bounded element magnitude, etc., on the system input matrix.

\end{abstract}

\IEEEpeerreviewmaketitle

\section{Introduction}

Controllability and observability are two fundamental properties required for a system. The former is concerned with possibilities of arbitrarily manoeuvring system internal states through external influences, while the latter with capabilities of estimating system internal states through measuring some system variables \cite{ksh,lsb11,Olshevsky14,zdg,zhou15}. These system properties have been extensively investigated for a long time, and various results have been obtained. There are, however, still many important issues needing further efforts. Among them, two essential problems are that to guarantee controllability/observability of a system, which states should be directly controlled/measured \cite{cd13,emccb12,wj01}. Due to the duality between controllability and observability of a system, these two problems are mathematically equivalent to each other.

With the development of network and communication technologies, etc., renewed interests have recently been widely attracted on large scale systems, which are now often called network systems \cite{emccb12,siljak78,sbkkmpr11}. Various long standing issues are being investigated from new perspectives, such as distributed estimation and control, analysis for their stability, controllability, observability, etc.\cite{hot11,np13,pka16,ra13,zhou13,zhou15,zz16}. On the other hand, many new important issues are being raised, such as structure identification, sparsity analysis, etc. \cite{csn09,Kolaczyk09,xz14}. Among these investigations, variable selections for these systems to insure controllability/observability are one of the most challenging issues\cite{emccb12,lsb11,pka16,scl15,zhou15}. In contrast to a traditional engineering system, a great amount of variables, for example, 1,000 or more, are usually involved in a network system. Moreover, there are great freedoms in selecting a set of variables to control or measure. Under such a situation, a brute force search, which  compares all possibilities through simply exhausting all variable combinations, are computationally prohibitive in the worst case, and systematic methods are required \cite{Olshevsky14,pka16,wj01,zhou15}.

Owing to various efforts of many researchers, some significant advancements have been recently achieved. In particular, it is proved in \cite{Olshevsky14} that some well known and widely investigated minimal controllability problems (MCP) are NP-hard, and a simple heuristic method generally works well for these problems that sequentially maximizes the rank increment in each step of the system controllability matrix. These problems are often met in actual engineering applications and the results are therefore of both significant theoretic values and great practical significance. Inspired by these results, various other results have been established, such as those on the problem of minimal input selection under average/worst case control energy restrictions reported, which are reported in \cite{scl15,pka16}, etc.

The key, and possibly also the most intelligent, step in \cite{Olshevsky14} is to include the well known hitting set problem (HSP) as a special MCP in which the state transition matrix has distinct eigenvalues and there is only one external input. Note that the HSP has been proved to be NP-hard, this embedment naturally leads to the NP-hardness of the MCP. In addition to these, that paper has also proved that three $k$-sparse controllability problems, including that associated with the previous MCP, are mathematically equivalent to each other, in which some have only one external input, while the others have several external inputs, also under the condition that all the eigenvalues of the state transition matrix are different from each other. Based on these results, it is declared that all the associated three MCPs are NP-hard.

These results have clearly settled some long standing important issues in control and system theories. The proofs about the problem equivalences, however, are based on a probabilistic point of view. In particular, a vector with some random variables has been generated to show the existence of a preferable input matrix, which can guarantee the existence of the required input matrix {\it only} with probability $1$. This further leads to that the equivalence of these problems can only be declared correct with probability one, also. Note that a random event with probability one is still not an event that must certainly occur. Moreover, in analyzing computational complexity of the HSP, a deterministic approach is adopted which reveals its computational costs in the worst case. These imply that a deterministic approach is more mathematically reasonable in these proofs, and it seems safe to declare that further efforts are still required to enhance the equivalence of these MCPs.

In this paper, we give another proof about the equivalence of these MCPs, which is equivalent to the equivalence of three associated $k$-sparse controllability problems. The derivations are based only on primary algebraic operations, and does not need any introduction of artificial random variables. These results enhance the results of \cite{Olshevsky14}, and appear more appropriate in computational complexity analysis for these MCPs, provided that it is combined with the HSP. In addition to these, this paper also makes it clear that even if there are some other restrictions on the system input matrix, such as that its elements are bounded in magnitude, etc., these MCPs are still mathematically equivalent to each other, which significantly extends applications of the results developed in \cite{Olshevsky14}.

The outline of this paper is as follows. In the next section, a description is given for the problem formulations. Section III develops main results of this paper, while Section V concludes this paper.

The following notation and symbols are adopted. ${\cal R}^{m\times n}$ and ${\cal C}^{m\times n}$ are utilized to represent the $m\times n$ dimensional real and complex linear spaces. When $m$ and/or $n$ are equal to $1$, they are usually omitted. $e_{k}$ stands for the $k$-th canonical basis vector of the linear space ${\cal R}^{n}$ or ${\cal C}^{n}$. ${\rm\bf diag}\!\{x_{i}|_{i=1}^{n}\}$ denotes a diagonal matrix with its $i$-th diagonal element being $x_{i}$, while ${\rm\bf col}\!\{x_{i}|_{i=1}^{n}\}$ the vector stacked by $x_{i}|_{i=1}^{n}$ with its $i$-th row element being $x_{i}$. ${\cal S}upp(x)$ is used to denote the support of the vector $x$, that is, the set consisting of all the row indices of this vector with its corresponding element not being zero. $\sharp({\cal A})$ represents the number of elements in the set $\cal A$, while $0_{m}$ and $0_{m\times n}$ respectively the $m$ dimensional column vector and the $m\times n$ dimensional matrix with all elements being zero. The superscript $T$ and $H$ are used to denote respectively the transpose and the conjugate transpose
of a matrix/vector, while $\bar{\cdot}$ the conjugate of a complex number/variable.

\section{Problem Formulation}

Consider three linear time invariant (LTI) dynamic systems ${\rm \bf \Sigma}_{i}$, $i=v,d,f$, with respectively the following state transition equations,
\begin{eqnarray}
& &\hspace*{-0.85cm} {\rm\bf\Sigma}_{v}:\hspace{0.1cm} \frac{{\rm d}x(t)}{{\rm d} t}=Ax(t)+B_{v}u(t)   \label{eqn:1} \\
& &\hspace*{-0.85cm} {\rm\bf\Sigma}_{d}:\hspace{0.1cm} \frac{{\rm d}x(t)}{{\rm d} t}=Ax(t)+B_{d}u(t)   \label{eqn:2} \\
& &\hspace*{-0.85cm} {\rm\bf\Sigma}_{f}:\hspace{0.1cm} \frac{{\rm d}x(t)}{{\rm d} t}=Ax(t)+B_{f}u(t)    \label{eqn:3}
\end{eqnarray}
in which both $x(t)$ and $B_{v}$ are $n$ dimensional real valued vectors, while $B_{d}$ a $n\times n$ dimensional diagonal real valued matrix, $B_{f}$ a $n\times p$ dimensional real valued matrix. The other vectors and matrices have compatible dimensions.

In system designs, an interesting problem is to find the sparest $B_{i}$, that is, a vector/matrix $B_{i}$ with the minimal number of nonzero elements, such that the corresponding LTI system ${\rm \bf \Sigma}_{i}$ is controllable, $i=v,d,f$. These problems are usually investigated under the condition that the integers $n$ and $p$ are fixed and the matrix $A$ is known, and are extensively called minimal controllability problems (MCP)\cite{Olshevsky14,pka16}

To avoid awkward statements, the following terminologies are adopted in this paper which are similar to those of \cite{Olshevsky14}. A vector/matrix is called $k$-sparse when it has at most $k$ nonzero elements. When a system is controllable with its state transition matrix and input matrix respectively being $A$ and $B$, it is sometimes also stated as the matrix pair $(A,B)$ is controllable. Similar statements are also adopted for the observability of a linear dynamic system.  If there exists a $k$-sparse vector/matrix $B_{i}$, such that the matrix pair $(A, B_{i})$ is controllable, then, we call the system ${\rm \bf \Sigma}_{i}$ $k$-sparse controllable, $i=v,d,f$.

In \cite{Olshevsky14}, it has been declared that when the matrix $A$ has distinct eigenvalues, these three MCPs are equivalently difficult in mathematics. The most important step there is to establish the equivalence of the $k$-sparse controllability of these three LTI dynamic systems. In its proofs, however, some random variables are introduced. This introduction makes the corresponding conclusions correct only in a probabilistic point of view, which may not be very appropriate in computational complexity for these MCPs.

This equivalence problem is re-investigated in this paper using a completely deterministic approach. That is, the purpose of this paper is to establish mathematical equivalence among these three MCPs, which is equivalent to the equivalence of the three associated $k$-sparse controllability problems, without introducing any random variables. In this investigation, we still assume that the state transition matrix $A$ in these systems does not have any repeated eigenvalues. This does not sacrifice any generality of the conclusions about the NP-hardness of these MCPs, as the HSP has been successfully embedded in \cite{Olshevsky14} into a MCP in which the system matrices always satisfy these assumption.

In this paper, we will also briefly discuss these MCPs under the situation that there are some other restrictions on the system input matrix $B_{i}$, $i=v,d,f$. For example, its element magnitude is bounded by a prescribed positive number, its Frobenius norm is not greater than some fixed value, etc., which is often met in actual engineering applications \cite{pka16,scl15,zhou15}.

\section{Main Results}

To investigate the equivalence of the aforementioned three MCPs, we at first need the following results on system controllability and observability, which is widely known in linear system theory, and extensively called as the PBH test \cite{ksh,zdg,zhou15}.

\hspace*{-0.45cm}{\bf Lemma 1.} Consider a continuous LTI system described by
\begin{displaymath}
\frac{{\rm d}x(t)}{{\rm d} t}=Ax(t)+Bu(t),\hspace{0.5cm}
y(t)=Cx(t)+Du(t)
\end{displaymath}
\begin{itemize}
\item The system is controllable if and only if for every complex scalar $\lambda$ and every nonzero complex vector $x$ satisfying $x^{H}A=\lambda x^{H}$, $x^{H}B\neq 0$.
\item The system is observable if and only if for every complex scalar $\lambda$ and every nonzero complex vector $y$ satisfying $Ay=\lambda y$, $Cy\neq 0$.
\end{itemize}

From this lemma, it is clear that the observability of the matrix pair $(A,\;C)$ is equivalent to the controllability of the matrix pair $(A^{T},\;B^{T})$. Therefore, results of this paper can be directly applied to a corresponding minimal observability problem, in which the sparest vector or matrix $C$ is searched under the restriction that the system is observable.

Note that for an arbitrary nonzero column vector, it is always possible to normalize its Euclidean norm into 1 while simultaneously to make its first nonzero element from the top positive. More precisely, assume that $x\neq 0$ and its $j$-th row element is the first nonzero element counting from the top. Then,
\begin{displaymath}
\frac{\bar{x}_{j}}{|x_{j}|\sqrt{x^{H}x}}x
\end{displaymath}
satisfies simultaneously these two conditions. Note also that a vector obtained through multiplying an eigenvector of a matrix by a constant is still an eigenvector of this matrix associated with the same eigenvalue. It can be declared that for each eigenvalue of a matrix, there exists at least one eigenvector, whose Euclidean norm is equal to one, while whose first nonzero element from the top is positive. In addition, it can be easily proved that for a square matrix with distinctive eigenvalues, associated with each of its eigenvalues, there is one and only one eigenvector that simultaneously holds these two properties.

On the basis of Lemma 1, the following conclusion is established, which gives a necessary and sufficient condition on the availability of an input matrix, such that the corresponding system is controllable.

\hspace*{-0.45cm}{\bf Theorem 1.} Assume that the $n\times n$ dimensional real matrix $A$ has $n$ distinct eigenvalues $\lambda_{i}$, $i=1,2,\cdots,n$. Let $x_{i}$ denote the left eigenvector associated with the eigenvalue $\lambda_{i}$, with its Euclidean norm equal to $1$ and its first nonzero element from the top positive. Then, for a prescribed set ${\cal S}_{v}$, which is a subset of the set $\{1,\;2,\;\cdots,\;n\}$, there exists a $n$ dimensional real vector $b$ with its support being a subset of the set ${\cal S}_{v}$, such that a continuous LTI system with its state transition equation being
\begin{equation}
\frac{{\rm d}x(t)}{{\rm d} t}=Ax(t)+bu(t)
\label{eqn:5}
\end{equation}
is controllable, if and only if for each $i\in \{1,\;2,\;\cdots,\;n\}$,
\begin{equation}
{\cal S}upp(x_{i})\bigcap {\cal S}_{v}\neq\emptyset \label{eqn:6}
\end{equation}

\hspace*{-0.45cm}{\bf\it Proof}:  Assume that there is a vector $b$ with ${\cal S}upp(b)={\cal S}_{v}$, such that a dynamic system with its state transition equation described by Equation (\ref{eqn:5}) is controllable. Moreover, assume that in the set $\{1,\;2,\;\cdots,\;n\}$,  there exists an integer $i$, such that
\begin{equation}
{\cal S}upp(x_{i})\bigcap {\cal S}_{v}=\emptyset
\end{equation}
Let $x_{i,j}$ and $b_{j}$ denote respectively the $j$-th row elements of the vectors $x_{i}$ and $b$. Then,
\begin{equation}
x_{i}^{H}b=\sum_{j=1}^{n}\bar{x}_{i,j}b_{j}=\sum_{j\in {\cal S}upp(x_{i})\bigcap {\cal S}_{v}}\bar{x}_{i,j}b_{j}=0
\end{equation}
This is clearly in contradiction with the results of Lemma 1. Therefore, the assumption about the existence of the aforemention integer $i$ is not reasonable.

On the contrary, assume that the condition of Equation (\ref{eqn:6}) is satisfied for every $1\leq i\leq n$, but there does not exist a vector $b$ with ${\cal S}upp(b)={\cal S}_{v}$, such that a dynamic system with its state transition equation described by Equation (\ref{eqn:5}) is controllable. For any prescribed $n$ dimensional real vector $b$, let ${\cal Z}_{b}$ denote the set consisting of the integers $i$ such that $x_{i}^{H}b=0$ with $1\leq i\leq n $. Mathematically speaking,
\begin{displaymath}
{\cal Z}_{b}=\left\{\; i\:\left|\:x_{i}^{H}b=0,\;i\in \{1\;2,\;\cdots,\;n\}\right.\right\}
\end{displaymath}
Using this set, define the set $\cal B$ as
\begin{equation}
{\cal B}=\left\{\: b\:\left|\: \sharp({\cal Z}_{b})\;\;{\rm is\;\; minimized},\; {\cal S}upp(b)={\cal S}_{v}\;\right.\right\}
\end{equation}
That is, this set consists of all the vectors $b$ with its support being ${\cal S}_{v}$ which minimizes the number of zero elements in the set $\{x_{i}^{H}b, i=1,2,\cdots, n\}$.

According to the assumptions, the eigenvalues of the matrix $A$, that is, $\lambda_{i}$, $i=1,2,\cdots, n$, are distinct. Then, the $n$ vectors ${x}_{i}$ are linearly independent of each other. Moreover, for any left eigenvector of this matrix, say $x$, there certainly exist an integer $i\in \{1,2 \cdots, n\}$ and a nonzero scalar $\alpha(i)$, such that
\begin{equation}
x^{H}A=\lambda_{i}{x}^{H}, \hspace{0.5cm} x=\alpha(i) x_{i}
\label{eqn:15}
\end{equation}
It can therefore be declared from Lemma 1 that when the eigenvalues of the matrix $A$ are distinct, the matrix pair $(A,\;b)$ is controllable, if and only if $x_{i}^{H}b\neq 0$ is valid for every $i=1,2,\cdots,n$.

Based on these observations, as well as the definitions of the set $\cal B$ and ${\cal Z}_{b}$, we have that for each $b\in {\cal B}$, it is certain that
\begin{equation}
1\leq \sharp({\cal Z}_{b}) \leq n
\end{equation}
that is, there exists at least one integer, denote it by $i$, such that $i\in \{1,2 \cdots, n\}$ and $x_{i}^{H}b=0$. Take any integer from the set
${\cal S}upp(b)\bigcap {\cal S}upp(x_{i})$, and denote it by $k(i)$. This operation is always possible, as from the assumption that the condition of Equation (\ref{eqn:6}) is satisfied for every $1\leq i\leq n$, we have that this intersection set is not empty.  Moreover, assume that in addition to the support of the vector $x_{i}$, $k(i)$ also belongs to the support of $x_{s_{j}(i,k)}$, in which $s_{j}(i,k)\in \{1,2 \cdots, n\}/\{i\}$ and $j=1,2,\cdots,\alpha(i,k)$. Denote $x_{s_{j}(i,k)}^{H}b$ by $\gamma_{j}(i,k)$, and define $\beta_{j}(i,k)$ as
\begin{displaymath}
\beta_{j}(i,k)=-\frac{\gamma_{j}(i,k)}{x_{s_{j}(i,k),\;k(i)}}
\end{displaymath}

Construct a set ${\cal P}(i,k)$ as
\begin{equation}
{\cal P}(i,k)={\cal R}\left/\left(\{0\}\bigcup\{\beta_{j}(i,k),\;j=1,2,\cdots, \alpha(i,k)\}\right)\right.
\label{eqn:16}
\end{equation}
Using an arbitrary $\delta_{b}$ which belongs to the set ${\cal P}(i,k)$, construct a $n$ dimensional vector $\hat{b}$ as
\begin{equation}
\hat{b}=b+\delta_{b}e_{k(i)}
\end{equation}
in which $e_{k(i)}$ is the $k(i)$-th canonical basis vector of the linear space ${\cal R}^{n}$.
Clearly, ${\cal S}upp(\hat{b})\subseteq {\cal S}upp(b)$. Hence, the support of the vector $\hat{b}$ is also included by the set ${\cal S}_{v}$. In addition,
\begin{eqnarray}
x_{i}^{H}\hat{b}&=&x_{i}^{H}\left(b+\delta_{b}e_{k(i)}\right) \nonumber\\
&=&x_{i,k(i)}\delta_{b} \nonumber\\
&\neq & 0
\end{eqnarray}
Moreover, for an arbitrary $s_{j}(i,k)$ with $j=1,2,\cdots,\alpha(i,k)$, we have that
\begin{eqnarray}
x_{s_{j}(i,k)}^{H}\hat{b}&=&x_{s_{j}(i,k)}^{H}\left(b+\delta_{b}e_{k(i)}\right) \nonumber\\
&=&\gamma_{j}(i,k)+x_{s_{j}(i,k),k(i)}\delta_{b} \nonumber\\
&\neq & 0
\end{eqnarray}
Furthermore, for each $m \in \{1,2,\cdots, n\}/\{i,s_{j}(i,k)|_{j=1}^{\alpha(i,k)}\}$, the following equality can be established from $k(i)\not\in {\cal S}upp(x_{m})$,
\begin{eqnarray}
x_{m}^{H}\hat{b}&=&x_{m}^{H}\left(b+\delta_{b}e_{k(i)}\right) \nonumber\\
&=&x_{m}^{H}b+x_{m,k(i)}\delta_{b} \nonumber\\
&=& x_{m}^{H}b
\end{eqnarray}

The above arguments remain valid if the integer $k(i)$ only belongs to the set ${\cal S}upp(x_{i})$, through simply modifying the set ${\cal P}(i,k)$ to be
${\cal P}(i,k)={\cal C}/\{0\}$.

Hence, by means of replacing the vector $b$ by the vector $\hat{b}$, the number of zero elements in the set ${\cal Z}_{b}$ can be reduced at least by one. More precisely
\begin{equation}
\sharp({\cal Z}_{\hat{b}})\leq  \sharp({\cal Z}_{b})-1
\end{equation}

This implies that $\sharp({\cal Z}_{b})$ does not achieve the minimal value by the selected vector $b$, which contradicts that this vector is chosen from the set consisting of all the vectors with its support being ${\cal S}_{v}$ which minimizes the number of zero elements in the set $\{x_{i}^{H}b, i=1,2,\cdots, n\}$. Hence, the assumption about the nonexistence of a vector $b$ such that the matrix pair $(A,b)$ is controllable, is not reasonable.

This completes the proof. \hspace{\fill}$\Diamond$

From the above derivations, it is clear that when the condition of Equation (\ref{eqn:6}) is satisfied, then, from each $n$ dimensional vector $b$ that has its support set equal to ${\cal S}_{v}$, it is possible to construct a vector of the same dimension, denote it by $\tilde{b}$, such as the matrix pair $(A,\tilde{b})$ is controllable. The reasons are that, for each $b\in{\cal R}^{n}$, $\sharp({\cal Z}_{b})$ takes a finite nonnegative value, which is in fact always not greater than $n$. More precisely, if $\sharp({\cal Z}_{b})$ happens to be zero, then, $\tilde{b}=b$ satisfies the requirements. On the other hand, if $\sharp({\cal Z}_{b})>0$, then, the procedure of replacing the vector $b$ by the vector $\hat{b}$ in the aforementioned proof on the sufficiency of the condition, can strictly reduces the number of zero elements in the set ${\cal Z}_{b}$. It can therefore be declared that by repeating this procedure at most $n$ times, a $n$ dimensional vector $\tilde{b}$ can be constructed, such that
\begin{equation}
\sharp({\cal Z}_{\tilde{b}})=0
\end{equation}
That is, the matrix pair $(A,\tilde{b})$ is controllable. However, it is worthwhile to mention that although the support of the vector $\tilde{b}$ is guaranteed to also be a subset of the set ${\cal S}_{v}$, this vector usually is not the sparest one.

In the above derivations, the assumption that the eigenvalues of the matrix $A$ are distinct is very essential. More precisely, if this assumption is not satisfied, the conclusions of Theorem 2 may no longer be valid, even when this matrix has $n$ linearly independent left eigenvectors $x_{i}|_{i=1}^{n}$. Note that if a matrix has more than one left eigenvectors that are linearly independent of each other and are associated with the same eigenvalue, then, each linear combination of these eigenvectors is also a left eigenvector of this matrix associated with the same eigenvalue \cite{ksh,zdg,zhou15}. This makes the 2nd equality of Equation (\ref{eqn:15}) sometimes not be satisfied by all the left eigenvectors of a square matrix with repeated eigenvalues, which further invalidates the results of Theorem 1. As a matter of fact, our investigations show that the condition of Equation (\ref{eqn:6}) is generally {\it only necessary} for the controllability of the matrix pair $(A,\;b)$.

Another thing worth of mentioning is that in \cite{Olshevsky14}, it has been clearly pointed out that "The argument is an application of the PBH test for controllability which tells us that to make (1) controllable, it is necessary and sufficient to choose a vector $b$ that is not orthogonal to all of the left-eigenvectors of the matrix $A$. It is easy to see that if $A$ does not have any repeated eigenvalues, this is possible if and
only if the support of $b$ intersects with the support of every left-eigenvector of $A$.". This is almost exactly the conclusion of Theorem 1, although it appears better to replace the word "all" with "each" in the expressions "$\cdots$ is not orthogonal to all of the left-eigenvectors of the matrix $A$". But there is no proof of these conclusions in \cite{Olshevsky14}. Interestingly, it will be seen soon that it is just the results of Theorem 1 that makes us possible to remove the artificial random variables adopted in the proof of that paper, which leads to arguments that are more appropriate in analyzing the computational complexity of the aforementioned MCPs.

Note that the support of a vector does not change after multiplying it with a nonzero constant. It is clear that when the state transition matrix $A$ does not have any  repeated eigenvalues, the support of each of its left eigenvectors depends neither on its Euclidean norm nor on the sign of its first nonzero element from the top. In order to avoid possible misunderstanding and awkward statements, however, the left eigenvectors $x_{i}|_{i=0}^{n}$ are required to satisfy these norm condition and sign condition, simply because that these eigenvectors are uniquely determined by the state transition matrix $A$.

It is interesting to note that obviously from the proof of Theorem 1, its results remains valid even if there are some other constraints on the system input matrix, such as magnitude restrictions on its elements and/or its Euclidean norm, etc., which are often met in actual applications \cite{pka16,scl15,zhou15}. For example, if each element of this matrix, which is in fact the vector $b$ in this situation, is required to be bounded in magnitude by a prescribed positive number, say, $h$, then, the proof of Theorem 1 remains valid through modifying {\it only} the set ${\cal P}(i,k)$ of Equation (\ref{eqn:16}) as
\begin{displaymath}
{\cal P}(i,k)=(-h,\; h)\left/\left(\{0\}\bigcup\{\beta_{j}(i,k),\;j=1,2,\cdots, \alpha(i,k)\}\right)\right.
\end{displaymath}
provided that the vector $b$ satisfies this restriction, which can be easily met when the set $\cal B$ is not empty, noting that if $b\in{\cal B}$, then $\alpha b$ also belongs to the set $\cal B$ for every positive number. Note that there are infinitely many elements in the set $(-h,\; h)$. This set will never be empty which is essential in searching the preferable vector $\tilde{b}$, as the procedure of replacing the vector $b$ by the vector $\hat{b}$ will be repeated at most $n$ times. When there are some restrictions on the Euclidean norm of the input matrix, the search of the preferable input matrix in the above sufficiency proof must be started from a vector satisfying these restrictions. Afterwards, these restrictions can be equivalently expressed as some constraints on the magnitude of the element $b_{k(i)}$, and the remaining derivations are completely the same as those when the input matrix is magnitude bounded element-wise. In fact, if the restrictions on the system input matrix does not make any of its elements being constrained to a set consisting of only some isolated values, then, it can be proved that the results of Theorem 1 remains valid. The details are omitted due to their obviousness.

Using these results, equivalences among the aforementioned three MCPs can be established. At first, we investigate relations between the MCPs with System ${\bf\Sigma}_{v}$ and System ${\bf\Sigma}_{d}$.

\hspace*{-0.45cm}{\bf Theorem 2.} Assume that the eigenvalues of the matrix $A$ are distinct. Then, the LTI dynamic system ${\rm \bf \Sigma}_{v}$ is $k$-sparse controllable, if and only if the LTI dynamic system ${\rm \bf \Sigma}_{d}$ is $k$-sparse controllable.

\hspace*{-0.45cm}{\bf\it Proof}: From the assumption that the eigenvalues of the matrix $A$ are distinct, it can be declared that this matrix has $n$ linearly independent eigenvectors. Let $x_{i}$, $i=1,2,\cdots,n$, denote its left eigenvectors which are linearly independent, whose Euclidean norms are equal to $1$, and whose  first nonzero elements from the top are positive.

Assume that the system ${\rm \bf \Sigma}_{v}$ is $k$-sparse controllable. Then, there exists at least one $k$-sparse vector $B_{v}$, that is, $\sharp({\cal S}upp(B_{v}))\leq k$, such that the matrix pair $(A, B_{v})$ is controllable. Then, according to Theorem 1, it is certain that for each $i=1,2,\cdots,n$,
\begin{equation}
{\cal S}upp(x_{i})\bigcap {\cal S}upp(B_{v})\neq\emptyset
\end{equation}
Equivalently, for an arbitrary $i\in \{1,2,\cdots,n\}$, there exists at least one $j(i)\in \{1,2,\cdots,n\}$, such that
\begin{equation}
x_{i,j(i)}\neq 0,\hspace{0.5cm} B_{v,j(i)}\neq 0
\label{eqn:7}
\end{equation}

Chose a arbitrary vector $b$ with $b={\rm\bf col}\!\{b_{i}|_{i=1}^{n}\}$ and ${\cal S}upp(b)={\cal S}upp(B_{v})$. Moreover, construct a diagonal matrix $B_{d}$ using this vector as $B_{d}={\rm\bf diag}\!\{b_{i}|_{i=1}^{n}\}$. Partition the matrix $B_{d}$ as
\begin{equation}
B_{d}=[b_{d,1}\;\; b_{d,2}\;\; \cdots\;\; b_{d,n}]
\end{equation}
in which $b_{d,i}$ is a $n$ dimensional column vector, $i=1,2,\cdots,n$.

For every $i\in\{1,2,\cdots,n\}$, we have from this construction and Equation (\ref{eqn:7}) that $b_{j(i)}\neq 0$, which further leads to
\begin{equation}
x_{i}^{H}b_{d,j(i)}=x_{i,j(i)}b_{d,j(i),j(i)}\neq 0
\end{equation}
Hence
\begin{equation}
x_{i}^{H}B_{d}\neq 0
\label{eqn:8}
\end{equation}

Recall that the eigenvalues of the matrix $A$ is different. We therefore have that for each of its left eigenvector, say, $x$, there certainly exists a $i\in\{1,2,\cdots,n\}$, and a nonzero complex scalar $\alpha(i)$, such that $x=\alpha(i) x_{i}$. It can therefore be declared from Equation (\ref{eqn:8}) that $x^{H}B_{d}\neq 0$. Hence, controllability of the matrix pair $(A,B_{d})$ can be claimed from Lemma 1. As ${\cal S}upp(B_{d})={\cal S}upp(b)={\cal S}upp(B_{v})\leq k$, we have that the system ${\rm \bf \Sigma}_{d}$ is also $k$-sparse controllable.

On the contrary, assume that the system ${\rm \bf \Sigma}_{d}$ is $k$-sparse controllable. Then, there exists a $k$-sparse diagonal matrix $B_{d}$, such that the matrix pair $(A,B_{d})$ is controllable. From Lemma 1 and the diagonal structure of the matrix $B_{d}$, it can be declared that for each $i\in\{1,2,\cdots,n\}$, there is at least one $j(i)\in\{1,2,\cdots,n\}$, such that
\begin{equation}
x_{i,j(i)}\neq 0,\hspace{0.5cm} B_{d,j(i),j(i)}\neq 0
\label{eqn:9}
\end{equation}
Define Sets ${\cal B}$ and ${\cal B}_{i}$, $i=1,2,\cdots,n$, respectively as
\begin{eqnarray}
& & \hspace*{-1cm} {\cal B}_{i}=\left\{\: k \: \left|\: x_{i,k}\neq 0,\hspace{0.15cm} B_{d,k,k}\neq 0, \hspace{0.15cm} k\in\{1,2,\cdots,n\} \:\right.\right\}
\label{eqn:10} \\
& & \hspace*{-1cm} {\cal B}=\bigcup_{i=1}^{n}{\cal B}_{i}
\label{eqn:11}
\end{eqnarray}
Then, from Equation (\ref{eqn:9}), it can be declared that the set ${\cal B}_{i}$ is not empty for each $i=1,2,\cdots, n$. Moreover, from the $k$-sparseness of the diagonal matrix $B_{d}$, it is obvious that $\sharp({\cal B})\leq k$.

Chose an arbitrary $n$ dimensional vector $b$ with its support equal to $\cal B$. Then, this vector is certainly $k$-sparse. Moreover, from the definitions of the sets ${\cal B}$ and ${\cal B}_{i}|_{i=1}^{n}$ given respectively in Equations (\ref{eqn:10}) and (\ref{eqn:11}),
we have that for each $i=1,2,\cdots,n$,
\begin{eqnarray}
& & {\cal S}upp(x_{i})\bigcap {\cal S}upp(b) \nonumber\\
&=&{\cal S}upp(x_{i})\bigcap {\cal B} \nonumber\\
&=& {\cal S}upp(x_{i})\bigcap \left(\bigcup_{j=1}^{n}{\cal B}_{j}\right)  \nonumber\\
&=& \bigcup_{j=1}^{n}\left({\cal S}upp(x_{i})\bigcap {\cal B}_{j}\right) \nonumber\\
&\supseteq & {\cal S}upp(x_{i})\bigcap {\cal B}_{i} \nonumber\\
&= & {\cal B}_{i} \nonumber\\
&\neq & \emptyset
\label{eqn:12}
\end{eqnarray}

It can therefore be declared from Theorem 1 that the system ${\rm \bf \Sigma}_{v}$ is also $k$-sparse controllable.

This completes the proof. \hspace{\fill}$\Diamond$

Note that in System ${\rm \bf \Sigma}_{v}$, there is {\it only} one actuator. But in System ${\rm \bf \Sigma}_{d}$, there are $n$ actuators. However, clearly from Theorem 2, when only controllability is concerned and the state transition matrix does not have any repeated eigenvalues, the number of actuators are not very important. The most important thing is the number of system states that can be directly manipulated. Note that it is physically intuitive and widely accepted that the more the actuators are, the easier the system is to be controlled \cite{lsb11}. This conclusion is a little surprising from an engineering point of view.

If the equivalence of Theorem 2 can be explained as a result of the equal freedoms in the matrices $B_{v}$ and $B_{d}$, both of them are in fact equal to that of selecting no more than $k$ positions from $n$ candidates, the following results show that even when the freedoms of input matrices are significantly different, minimal controllability of a system may still remains unchanged. These results may imply that while controllability is an important system property, an appropriate metric is needed for quantitatively evaluating how difficult a system is to be controlled. Similar observations have also been found in \cite{pka16,scl15,zhou15}.

\hspace*{-0.45cm}{\bf Theorem 3.}The LTI dynamic system ${\rm \bf \Sigma}_{v}$ is $k$-sparse controllable, if and only if the LTI dynamic system ${\rm \bf \Sigma}_{f}$ is $k$-sparse controllable, provided that the matrix $A$ has different eigenvalues.

\hspace*{-0.45cm}{\bf\it Proof}:  Assume that the system ${\rm \bf \Sigma}_{v}$ is $k$-sparse controllable. Then, there exists a $k$-sparse vector $B_{v}$ such that the matrix pair $(A,B_{v})$ is controllable. Hence, according to Lemma 1, for any nonzero vector $x$ satisfying $x^{H}A=\lambda x^{H}$ in which $\lambda$ is a complex scalar, we have that $x^{H}B_{v}\neq 0$. Define a matrix $B_{f}$ as
\begin{equation}
B_{f}=[0_{n\times (p-1)}\;\; B_{v}]
\end{equation}
Then, obviously $B_{f}\in {\cal R}^{n\times p}$ and is $k$-sparse. Moreover, note that for an arbitrary complex vector with a compatible dimension, say, $y$, we have that
\begin{equation}
y^{H}B_{f}=[0_{1\times (p-1)}\;\; y^{H}B_{v}]
\end{equation}
It can therefore be declared from the controllability of the matrix pair $(A,B_{v})$ and Lemma 1, that the matrix pair $(A,B_{f})$ is also controllable. Hence,  the system ${\rm \bf \Sigma}_{f}$ is $k$-sparse controllable.

On the contrary, assume that the system ${\rm \bf \Sigma}_{f}$ is $k$-sparse controllable. Then, there exists at least one $n\times p$ dimensional matrix $B_{f}$ with less than  $k$ nonzero elements, such that the matrix pair $(A,B_{f})$ is controllable. Partition the matrix $B_{f}$ as
\begin{displaymath}
B_{f}=\left[b_{f,1}\;\; b_{f,2}\;\; \cdots \;\; b_{f,p}\right]
\end{displaymath}
in which $b_{f,j}$ is a $n$ dimensional real vector, $j=1,2,\cdots,p$. Then, from the assumption on the matrix $B_{f}$, it is obvious that
\begin{equation}
\sum_{i=1}^{p}{\cal S}upp(b_{f,i})\leq k
\label{eqn:13}
\end{equation}
Moreover, according to Lemma 1 and the controllability of the matrix pair $(A,B_{f})$, we have that for each $i\in\{1,\;2,\;\cdots,\;n\}$, there exists at least one
$j(i)\in \{1,\;2,\;\cdots,\;p\}$, such that
\begin{equation}
x_{i}^{H}b_{f,j(i)}\neq 0
\label{eqn:14}
\end{equation}
As in Theorem 1, $x_{i}$, $i=1,2,\cdots,n$, here again denote the $n$ linearly independent left eigenvectors of the matrix $A$, whose Euclidean norms equal to $1$ and whose first nonzero elements from the top are positive.

For each $i=1,2,\cdots,n$, define a set ${\cal J}(i)$ as
\begin{equation}
{\cal J}(i)=\left\{\:k \: \left|\: x_{i}^{H}b_{f,k}\neq 0,\; 1\leq k\leq p\:\right.\right\}
\end{equation}
Then, it can be claimed from Equation (\ref{eqn:14}) that the set ${\cal J}(i)$ is not empty whenever $i\in\{1,\;2,\;\cdots,\; n\}$.

Using these definitions, define another set $\cal B$ as
\begin{equation}
{\cal B}=\bigcup_{i=1}^{n}\bigcup_{j\in{\cal J}(i)}{\cal S}upp(b_{f,j})
\end{equation}
Then, from Equation (\ref{eqn:13}), we have that for any $n$ dimensional vector $b$ with ${\cal S}upp(b)={\cal B}$,
\begin{eqnarray}
\sharp({\cal S}upp(b))&=& \sharp({\cal B})  \nonumber\\
&\leq & \sum_{i=1}^{p}{\cal S}upp(b_{f,i})  \nonumber\\
&\leq & k
\end{eqnarray}
Moreover, from the definition of the set ${\cal J}(i)$ and Equation (\ref{eqn:14}), it can be declared that for each $i=1,2,\cdots,n$,
\begin{eqnarray}
& & {\cal S}upp(x_{i}) \bigcap {\cal S}upp(b) \nonumber\\
&=& {\cal S}upp(x_{i}) \bigcap \left\{\bigcup_{q=1}^{n}\bigcup_{j\in{\cal J}(q)}{\cal S}upp(b_{f,j})\right\} \nonumber\\
&=& \bigcup_{q=1}^{n}\bigcup_{j\in{\cal J}(q)}\left\{ {\cal S}upp(x_{i}) \bigcap {\cal S}upp(b_{f,j})\right\} \nonumber\\
&\supseteq & \bigcup_{j\in{\cal J}(i)}\left\{ {\cal S}upp(x_{i}) \bigcap {\cal S}upp(b_{f,j})\right\} \nonumber\\
&\neq & \emptyset
\end{eqnarray}

On the basis of these results and Theorem 1, the existence of a $k$-sparse vector $b$ can be claimed that makes the matrix pair $(A,b)$ controllable. Hence, the system ${\rm \bf \Sigma}_{v}$ is $k$-sparse controllable.

This completes the proof. \hspace{\fill}$\Diamond$

It is worthwhile to point out here that although the expressions are different, the basic ideas in the first parts of the proofs of Theorems 2 and 3, that is, the "necessity" parts, are similar to those of \cite{Olshevsky14}. They are included in this paper for the completeness of the proof. The second parts of these proofs, that is, the "sufficiency" parts, however, are completely different. Here, the derivations are established on Theorem 1, which has been briefly mentioned in \cite{Olshevsky14}, but has not been rigorously proved there. The most important thing here is that the approach taken here is completely deterministic, while some artificial random variables are utilized in \cite{Olshevsky14} that makes the conclusions there valid only in a probabilistic point of view.

More precisely, when the system ${\bf \Sigma}_{d}$ or the system ${\bf \Sigma}_{f}$ is $k$-sparse controllable, in order to guarantee the existence of a $k$-sparse vector $b$ such that the matrix pair $(A,B_{v})$ is controllable, a random vector with elements being normally distributed is constructed in \cite{Olshevsky14}. This makes the $k$-sparse controllability of the system ${\bf \Sigma}_{v}$ be assured {\it only} with probability $1$, which does not match very well with the methodologies adopted in the NP-hardness of the HST problem. It appears therefore safe to declare that the approach taken in this paper is mathematically much easier to be understood, and is more suitable for computational complexity analysis of the associated optimization problem.

The following results can be immediately obtained from Theorems 2 and 3.

\hspace*{-0.45cm}{\bf Corollary 1.} Assume that all the eigenvalues of the matrix $A$ are different. Then, the LTI dynamic system ${\rm \bf \Sigma}_{d}$ is $k$-sparse controllable, if and only if the LTI dynamic system ${\rm \bf \Sigma}_{f}$ is $k$-sparse controllable,

\hspace*{-0.45cm}{\bf\it Proof}:  When the state transition matrix $A$ has distinctive eigenvalues, it has been proved respectively in Theorems 2 and 3 that, the $k$-sparse controllability of the system ${\rm \bf \Sigma}_{v}$ is equivalent both to that of the system ${\rm \bf \Sigma}_{d}$ and to that of the system ${\rm \bf \Sigma}_{f}$. Hence, the $k$-sparse controllability of the system ${\rm \bf \Sigma}_{d}$ and that of the system ${\rm \bf \Sigma}_{f}$ are also equivalent to each other. This completes the proof. \hspace{\fill}$\Diamond$

As being emphasized after the proof of Theorem 1, the assumption plays an essential role in all the above derivations that the state transition matrix $A$ has distinctive eigenvalues. When this assumption is not satisfied, our preliminary investigations show that the conclusions might be significantly different. As the derivations are tedious and the results are not relevant to investigations on the computational complexity of these three MCPs, they are not included here. However, as the hitting set problem has been successfully embedded into a MCP in \cite{Olshevsky14} with its state transition matrix having distinctive eigenvalues,
the conclusions on the NP-hardness of the aforementioned three MCPs remain valid even when the matrix $A$ has some eigenvalues that are equal to each other. In fact, this embedment is the most essential idea given in \cite{Olshevsky14} in the proof of the NP-hardness of a MCP, and can be regarded to be highly intelligent.

As having pointed before, the results of Theorem 1 are also valid when there are some other constraints on the system input matrix, such as restrictions on its element magnitude and/or Frobenius norm, etc. From the proofs of Theorems 2 and 3, as well as that of Corollary 1, it is clear that when constraints like these are also put on the input matrix $B_{i}$ of the system ${\bf\Sigma}_{i}$, $i=v,d,f$, the equivalence of the discussed three MCPs is still valid. The only requirement on these constraints are that each of the element of the input matrix should not be restricted to a set with only some isolated values. The details are omitted for avoiding awkward statements.

\section{Concluding Remarks}

In this paper, we have re-investigated the equivalence of three minimal controllability problems discussed in \cite{Olshevsky14}, using a completely deterministic approach. The equivalence has been established again without introducing any random variables. In addition, it has also been made clear that this equivalence remains valid even if there are some other constraints on the system input matrix, provided that these constraints do not make any of its elements belongs to a set consisting of only some isolated values. The results can be directly transformed to corresponding minimal observability problems. Compared to the available derivations given in \cite{Olshevsky14}, the derivations reported in this paper are mathematically much easier to be understood and appear more appropriate in analyzing computational complexity of the associated minimal controllability problems. When combined with Theorem 1 of \cite{Olshevsky14}, these results can declare the NP-hardness of these minimal controllability problems, as well as their approximations.

\end{document}